\documentclass[11pt,reqno,a4paper]{amsart}
\usepackage[margin=1in]{geometry}
\usepackage[usenames]{color}
\usepackage{amsmath,pdfsync,verbatim,graphicx,epstopdf,enumerate}
\usepackage{placeins}
\usepackage[abbrev,nobysame]{amsrefs}
\usepackage{float}
\usepackage[colorlinks=true]{hyperref}
\usepackage{cancel}
\usepackage[framemethod=tikz]{mdframed}
\hypersetup{allcolors=blue}
\usepackage{subcaption}
\usepackage{graphicx}
\allowdisplaybreaks
\numberwithin{equation}{section}

\newcommand{\lb}{\left(}

\newcommand{\rb}{\right)}

\newcommand{\Beq}{\begin{equation}}
	\newcommand{\Eeq}{\end{equation}}
\newcommand{\beq}{\begin{equation*}}
	\newcommand{\eeq}{\end{equation*}}
\newcommand{\bal}{\begin{align}}
	\newcommand{\eal}{\end{align}}

\usepackage{mathtools}
\mathtoolsset{showonlyrefs}

\usepackage[notref,notcite]{}

\newcommand{\bp}{\begin{prob}}
	\newcommand{\ep}{\end{prob}}
\newcommand{\bpr}{\begin{proof}}
	\newcommand{\epr}{\end{proof}}

\newcommand{\bel}[1]{\begin{equation}\label{#1}}
	\newcommand{\ee}{\end{equation}}

\newtheorem{theorem}{Theorem}[section]

\newtheorem{lemma}[theorem]{Lemma}

\newtheorem{example}[theorem]{Example}

\theoremstyle{definition}

\newtheorem{remark}[theorem]{Remark}

\newcommand{\D}{\mathrm{d}}

\usepackage[usenames]{color}
\usepackage{amsmath,pdfsync,verbatim,graphicx,epstopdf,enumerate}

\newcommand{\Rc}{\mathcal{R}}

\newcommand{\Sb}{\mathbb{S}}

\begin{document}
\title[Explicit inversion of SRT]{Explicit inversion of spherical Radon transforms in odd dimensions with partial radial data}
\author[Chatterjee, Krishnan and Tushir]{Pradipta Chatterjee, Venkateswaran P.\ Krishnan and Abhilash Tushir}
\address{
	Tata Institute of Fundamental Research, Centre For Applicable Mathematics
    \endgraf Bangalore, Karnataka 560065, India
	\endgraf
    \it Email: \tt pradipta22@tifrbng.res.in, vkrishnan@tifrbng.res.in, abhilash2296@gmail.com 
}
 \begin{abstract}
We derive an explicit inversion algorithm for the spherical Radon transform in odd dimensions with partial radial data. We prove that the reconstruction of the unknown function can be reduced to solving ordinary differential equations, thereby providing a more explicit approach in odd dimensions than solving Volterra integral equation of the first kind established in prior works. We also provide analytical  solutions in some special cases. Finally, we present numerical simulations  validating our theoretical results. Our work answers a question posed by Rubin in ``Inversion formulae for the spherical mean in odd dimensions and the Euler-Poisson-Darboux equation,''
Inverse Problems 24 (2008), no. 2, 025021, 10 pp. 
 \end{abstract}
\maketitle 
 \section{Introduction}
The spherical mean transform (SMT) of a function $f$ 
is defined as
\begin{equation}\label{SRTdef}
    \mathcal{R} f(p,t)=\frac{1}{\omega_{n-1}}\int\limits_{\mathbb{S}^{n-1}} f(p+t\theta) \mathrm{d} S (\theta), 
\end{equation}
where $\omega_{n-1}$ and $\mathrm{d}S$ represent the surface area and surface measure of $\mathbb{S}^{n-1}$, respectively. 

The study of SMT has attracted considerable attention due to its theoretical importance as well as in applications.
For example, SMT appears in the analysis of the wave equation, Euler-Poisson-Darboux equation and related PDEs; see  \cite{CH_Book,John-book,Rhee}.  Furthermore, SMT has been examined in the context of problems that appear in integral geometry and approximation theory; see \cite{agranovsky1996approximation,agranovsky1996injectivity,AQ2, AQ3, Agranovsky:1999,And,aramyan2020recovering} and the references therein. 

In applications, SMT appears in the mathematical study of photoacoustic ultrasound reconstruction, thermoacoustic tomography,  seismic imaging, quantum and optical imaging etc to name a few; see \cite{Kruger:1995,Kruger:V2,xu2002time,Xu:2005}. The basic principle of these reconstruction approaches involves illuminating the object with a short burst of radio-frequency waves, resulting in rapid, but minor, thermal expansion that produces a pressure wave. The pressure wave can be analyzed at the boundary or exterior of the object and is then utilized to reconstruct the function.

Several results concerning the inversion of SMT exist under different constraints on the centers, radii, and available data. A majority of the inversion formulas employs the spherical acquisition geometry, that is, where the centers of the integration spheres lie on a sphere encompassing the object of interest. In this setting, series-type inversion formulas were established in two-dimensions and three-dimensions by Norton and Norton, Linzer, Xu and Wang in \cite{norton1980reconstruction,norton1981ultrasonic,xu2002time}. In \cite{Finch-P-R}, Finch, Patch and  Rakesh derived backprojection-type inversion formulas in odd dimensions and later Finch, Haltmeier and Rakesh established similar inversion formulas for the case of even dimensions in \cite{Finch-Haltmeir-Rakesh_even-inversion}. Rubin in \cite{R} gave a straightforward proof of the Finch–Patch–Rakesh inversion formula in odd dimensions using the argument of analytic continuation and Erdelyi–Kober fractional integrals.  We also refer a work of Kunyansky \cite{Kunyansky}, where he derived an inversion for SMT by analyzing the Helmholtz equation. 

The majority of the aforementioned papers deal with full data in the radial direction. However, full radial data may not be available in some imaging applications such as imaging the area surrounding a bone. 
The primary purpose of this paper is to address a question posed by Rubin in \cite{R}. The question specifically asks for a reconstruction algorithm for SMT in odd dimensions in the aforementioned geometry and with partial radial data, that is, where the SMT data is available only for all $t$ with  $0<t<1$. This question was independently investigated, first by Ambartsoumian, Gouia-Zarrad and Lewis for SMT in 2-dimensions using Cormack's integral equations approach; see \cite{Ambartsoumian-Zarrad-Lewis}.  Using the same approach there has been several related works \cite{Ambartsoumian2015, Ambartsoumian2018,Salman_Article}. Numerical approaches of these algorithms have been studied as well \cite{Souvik:2015, Ambartsoumian2018}. The technique in all these papers has been to transform the inversion of SMT with partial radial data to a Volterra integral equation of the first kind with possibly a weakly singular kernel. This is then transformed to a Volterra integral equation of the second kind which is inverted using Picard's process of successive approximations, making the inversion algorithm non-explicit. 

We revisit the question of inversion of SMT with partial radial data in this work. Our main contribution is that in odd dimensions, the inversion of SMT with partial radial data can be made substantially more explicit and simpler compared to  prior works; it involves solving linear ODEs, thereby giving an affirmative answer to the question posed by Rubin \cite{R}. While our approach still lacks an explicit inversion formula (see \cite{Finch-P-R,R, Kunyansky} when full radial and angular data  of SMT are known), for the case of partial radial data in odd dimensions, to the best of our knowledge,  the current work seems to be the best available. In fact, in special cases, the solution to these ODEs can be explicitly found as well.  We deal with a few of these cases in this paper. Moreover, the inversion algorithm is computationally simpler, which makes it feasible and straightforward to apply in numerical reconstructions. We give a few reconstructions in this paper. 

The rest of the paper is organized as follows. 
Section \ref{sec:main result} provides the main results of this article. Section \ref{sec: preparatory results} gives a few lemmas required in the proofs of the main result.    Section \ref{sec:proof of main results} presents the proofs of main results.
In Section \ref{sec:analytical sol}, we provide the inversion formulas along with analytical solutions in certain special cases. 
Section \ref{sec:num simulations} presents a number of numerical simulations that illustrate the validity of the proposed inversion algorithm. 
\section{Main results}\label{sec:main result} 
 In this section, we give the main results proved in this article. 
 Our first result is an explicit inversion algorithm for SMT of radial functions. 
\begin{theorem}[Explicit inversion for radial functions]\label{mainthm1} Let $\epsilon>0$, $n\geq 3$ be an odd integer, and $k:=\frac{n-3}{2}$. Also assume that $f\in C_{c}^{\infty}(\mathbb{B})$ and  $$h_{k}(t)=\dfrac{4^{k}\omega_{2k+2}t^{2k+1}}{\omega_{2k+1}}\mathcal{R}f(p,t),~~\forall~(p,t)\in \mathbb{S}^{n-1}\times (0,1),$$   then $f$ can be recovered in the annular region $\mathbb{B}(\epsilon,1)$  by solving
\begin{multline*}
			\left[\frac{\D}{\D t}D^{2k}h_{k}\right](t)\\
			=(-1)^{k}k!4^{k}\sum_{m=0}^{k}\sum_{n=m}^{k}\sum_{l=n}^{k}2^{l+m-n}\frac{(2k-l)!}{(k-l)!m!}\binom{l+1}{n+1}\binom{l+n-m}{n-m}\frac{1}{t^{l}}\frac{(1-t)^{n+1}}{t^{n-m}}f^{(m)}(1-t),
		\end{multline*}
        with initial conditions $f^{(i)}(1-\epsilon^{\prime})=0$ for $0\leq i\leq k-1$, where $\epsilon^{\prime}>0$ is small enough such that $1-\epsilon^{\prime}$ lies outside the support of $f$.
\end{theorem}
\begin{remark}\label{remark:First}
    For $k=0$, the above theorem says that $f$ can be recovered by solving $\frac{\D}{\D t}h_{0}(t)=(1-t)f(1-t)$, and since no ODE is present, there is no requirement of any initial condition for $k=0$. This is given in \cite{R} as well.
\end{remark}

We now consider the general case. 
The spherical harmonics expansion of $f\in C_{c}^{\infty}(\mathbb{B})$ is given by
\begin{equation*}
        f(x)=\sum_{q=0}^{\infty}\sum_{s=1}^{d_q} f_{q,s}(|x|)Y_{q,s}\left(\frac{x}{|x|}\right),
\end{equation*}
where
\begin{equation}\label{dqandfqs}
 \quad d_{0}=1,\quad d_{q}=\frac{(2q+n-2)(n+q-3)}{q!(n-2)!},\quad \text{and}\quad f_{q,s}(r)=\int_{\mathbb{S}^{n-1}}f(r\theta)\overline{Y}_{q,s}(\theta)\D \theta.
\end{equation}
Further, the spherical harmonics expansion of $g=\mathcal{R}f$ is given by    
\begin{equation*}
     g(\theta,t)=\sum_{q=0}^{\infty}\sum_{s=1}^{d_q} g_{q,s}(t)Y_{q,s}(\theta),
\end{equation*}
where $g_{q,s}\in C_{c}^{\infty}((0,2))$.
Our main theorem for the explicit inversion of SMT for general functions is given as follows.
\begin{theorem}[Explicit inversion for general functions]\label{mainthm2}  Let $\epsilon>0$, $n\geq 3$ be an odd integer, and $k:=\frac{n-3}{2}$. Also assume that $f\in C_{c}^{\infty}(\mathbb{B})$ and for each $(q,s)$ satisfying $q\geq 0$ and $1\leq s\leq d_{q}$, we have
$$h_{q,s}(t):=\dfrac{\omega_{2k+2}t^{2k+1}}{\omega_{2k+1}}g_{q, s}(t)\quad \text{and}\quad \widetilde{f}_{q,s}(t)=\frac{f_{q,s}(t)}{t^q},$$
then $\widetilde{f}_{q,s}(t)$ can be recovered in annulus region $\mathbb{B}(\epsilon,1)$ by solving    \begin{multline*}
	\left[\frac{\D}{\D t}D^{q+2k}h_{q,s}\right](t)\\
=\frac{(-1)^{q+k}k!}{2^{q}}\sum_{m=0}^{q+k}\sum_{n=m}^{q+k}\sum_{l=n}^{q+k}2^{l+m-n}\frac{(2q+2k-l)!}{(q+k-l)!m!}\binom{l+1}{n+1}\binom{l+n-m}{n-m}\frac{1}{t^{l}}\frac{(1-t)^{n+1}}{t^{n-m}}\widetilde{f}_{q,s}^{(m)}(1-t),
\end{multline*}
with initial conditions $\widetilde{f}_{q,s}^{(i)}(1-\epsilon^{\prime})=0$ for $0\leq i\leq q+ k-1$, where $\epsilon^{\prime}>0$ is small enough such that $1-\epsilon^{\prime}$ lies outside the support of $f$.
\end{theorem}

\section{Preparatory results}\label{sec: preparatory results}
We provide the proof of the main theorems in the next section after developing some preparatory results to set the background.

For  radial functions, from \eqref{SRTdef},
\[
\Rc f(p,t)=\frac{1}{\omega_{n-1}}\int\limits_{\Sb^{n-1}} f(|p+t\theta|) \D S(\theta).
\]
Expanding this,
\[
\Rc f(p,t)=\frac{1}{\omega_{n-1}}\int\limits_{\Sb^{n-1}} f(\sqrt{1+ t^2 +2 t p\cdot \theta}) \D S(\theta).
\]
Let us apply Funk-Hecke theorem \cite{Seeley}. We bear in mind is that the support condition of $f$ forces restriction on the values that $p\cdot \theta$ can take. Note that $(-p) \cdot \theta$ ranges between $t/2$ and $1$. 
Hence 
\[
\Rc f(p,t)=\frac{\omega_{n-2}}{\omega_{n-1}}\int\limits_{t/2}^{1} f(\sqrt{1+ t^2 -2 st})(1-s^{2})^{\frac{n-3}{2}} \D s.
\]
With the change of variable $u=\sqrt{1+t^2-2st}$, we get
\begin{align*}
\Rc f(p,t) &= \frac{\omega_{n-2}}{\omega_{n-1}t}\int\limits_{1-t}^{1} uf(u) \left(1-\left(\frac{1+t^2-u^{2}}{2t}\right)^{2}\right)^{\frac{n-3}{2}} \D u\nonumber\\
&=\frac{\omega_{n-2}}{\omega_{n-1}2^{n-3}t^{n-2}}\int\limits_{1-t}^{1} uf(u) \left(4t^{2}-(1+t^2-u^{2})^{2}\right)^{\frac{n-3}{2}} \D u\nonumber\\
&=\frac{C_{n}}{t^{n-2}}\int\limits_{1-t}^{1} uf(u) [Q(t,u)]^{\frac{n-3}{2}} \D u,
\end{align*}
where
$$C_{n}=\frac{\omega_{n-2}}{\omega_{n-1}2^{n-3}}\quad\text{and}\quad Q(t,u)=((1+t)^2-u^2)(u^2-(1-t)^2).$$ 
Denoting $h(t)=\dfrac{t^{n-2}}{C_{n}}\mathcal{R}f(t)$ (note that since $f$ is radial, $\mathcal{R}f$ is independent of $p$), we get
\begin{equation*}
    h(t)=\int\limits_{1-t}^{1} uf(u) [Q(t,u)]^{\frac{n-3}{2}} \D u.
\end{equation*}
If we set $k:=\dfrac{n-3}{2}$ and using the notation $h_{k}(t)$ instead of $h(t)$ to emphasize  the dependence on $k$, then we have
\begin{equation}\label{hkdef}
    h_{k}(t)=\int\limits_{1-t}^{1} uf(u) [Q(t,u)]^{k} \D u.
\end{equation}

Our goal is to recover the function $f$ from the above expression. Since we are considering SMT in odd dimensions, $k$ is an integer and  we can differentiate the aforementioned expression several times to accomplish this task. We now gather certain information regarding the derivatives of $h_k$. 
We denote 
\begin{equation*}
    G_{i,j}(t):=\int\limits_{1-t}^{1}uf(u)[Q(t,u)]^{i}(u^2+1-t^2)^{j}\D u,\quad (i,j)\in \mathbb{N}_{0}^{2}.
\end{equation*}

Let us summarize a few immediate observations related to the derivatives of $G_{i,j}(t)$.
By using the fact that $Q(t,u)$ vanishes at $1-t$, the following conclusions can be drawn instantly:
\begin{enumerate}
\item $G_{k,0}(t)=h_{k}(t)$;\medskip
     \item $\frac{\D}{\D t}G_{0,0}(t)=(1-t)f(1-t)$;\medskip
		\item\label{rec:Gk0} $[DG_{i,0}](t)=4iG_{i-1,1}(t)$, for $i\geq 1$; \medskip
		\item\label{rec:Gkj} $[DG_{i,j}](t)=4iG_{i-1,j+1}(t)-2jG_{i,j-1}(t)$ for $i,j\geq 1$; and \medskip
        \item\label{rec:G0j} $[DG_{0,j}](t)=2^{j}T_{j}(t)-2jG_{0,j-1}(t)$ for $j\geq 1$, where $T_{j}(t)=\dfrac{(1-t)^{j+1}}{t}f(1-t)$.
\end{enumerate} 

\medskip 

The following lemma gives the formula for the  $(j+1)^{\mathrm{th}}$ derivative of $G_{0,j}(t)$.
\begin{lemma}\label{recurssion} For $j\geq 1$, the following relation holds:
\begin{equation}\label{rec:dDjGoj}
    \left[\frac{\D}{\D t}D^{j}G_{0,j}\right](t)
	=2^{j}\left((-1)^{j}j!(1-t)f(1-t)+\sum_{l=1}^{j}(-1)^{j-l}\frac{j!}{l!}\left[\frac{\D}{\D t}D^{l-1}T_{l}\right](t)\right),
\end{equation}
where $$T_{l}(t)=\dfrac{(1-t)^{l+1}}{t}f(1-t).$$
\end{lemma}
\begin{proof}
Let us prove using induction. For $j=1$, from \eqref{rec:G0j} we have
\begin{equation}
      [DG_{0,1}](t)
    =2\left(-G_{0,0}(t)+T_{1}(t)\right).
\end{equation}
Taking a $\frac{\D}{\D t}$ derivative, we immediately obtain \eqref{rec:dDjGoj} for $j=1$. Assuming that the relation \eqref{rec:dDjGoj} holds for $j$ and let us prove for $j+1$.  From \eqref{rec:G0j}, we have
 \begin{equation*}
 	\left[DG_{0,j+1}\right](t)
 	=2^{j+1}T_{j+1}(t)-2(j+1)G_{0,j}(t).
 \end{equation*}
 Now taking
$j$ $D$-derivatives and a $\frac{\D}{\D t}$ derivative, we get
\begin{equation*}
    \left[\frac{\D}{\D t}D^{j+1}G_{0,j+1}\right](t)=2^{j+1}\left[\frac{\D}{\D t}D^{j}T_{j+1}\right](t)-2(j+1)\left[\frac{\D}{\D t}D^{j}G_{0,j}\right](t).
\end{equation*}
 Using the induction hypothesis, the above expression takes the form
\begin{align*}
	&\left[\frac{\D}{\D t}D^{j+1}G_{0,j+1}\right](t)\\
&=2^{j+1}\left[\frac{\D}{\D t}D^{j}T_{j+1}\right](t)-2(j+1)2^{j}\left((-1)^{j}j!(1-t)f(1-t)+\sum_{l=1}^{j}(-1)^{j-l}\frac{j!}{l!}\left[\frac{\D}{\D t}D^{l-1}T_{l}\right](t)\right)\\
	&=2^{j+1}\left[\frac{\D}{\D t}[D^{j}T_{j+1}\right](t)\\&+2^{j+1}\left((-1)^{j+1}(j+1)!(1-t)f(1-t)+\sum_{l=1}^{j}(-1)^{j+1-l}\frac{(j+1)!}{l!}\left[\frac{\D}{\D t}D^{l-1}T_{l}\right](t)\right)\\
    &=2^{j+1}\left((-1)^{j+1}(j+1)!(1-t)f(1-t)+\sum_{l=1}^{j+1}(-1)^{j+1-l}\frac{(j+1)!}{l!}\left[\frac{\D}{\D t}D^{l-1}T_{l}\right](t)\right).
\end{align*}
Thus \eqref{rec:dDjGoj} holds for $j+1$. This completes the proof.
\end{proof}
Our goal is now to compute $[D^{l-1}T_{l}](t)$, but before we can apply the Leibniz rule to $T_{l}(t)=\dfrac{(1-t)^{l+1}}{t}f(1-t)$, we give the formulas for higher order $D$-derivatives of $\dfrac{(1-t)^{l+1}}{t}$ and that of $f(1-t)$ in terms of ordinary derivatives. Using these in \eqref{rec:dDjGoj}, we obtain the exact expression for $\left[\frac{\D}{\D t}D^{j}G_{0,j}\right](t)$ in terms of ordinary derivatives of $f$. 
\begin{lemma}\label{Tlprelim}
 For $p,l\geq 0$, we have
\Beq \label{tlder}
	\begin{aligned}
		D^{l}\left\{\frac{(1-t)^{p+1}}{t}\right\}&=(-1)^{l}\sum\limits_{s=0}^{l}\frac{(l+s)!(p+1)!}{2^{s}s!(l-s)!(p-l+s+1)!}\frac{(1-t)^{p-l+s+1}}{t^{l+s+1}}\\
        &=(-1)^{l}l!\sum\limits_{s=0}^{l} \frac{1}{2^s}{p+1\choose l-s}{l+s\choose s}\frac{(1-t)^{p-l+s+1}}{t^{l+s+1}}. 
	\end{aligned}
    \Eeq 
\end{lemma}
The proof is a straightforward verification. We skip it. 
\begin{lemma}\label{fder} For any $r\in\mathbb{N}$ and $f\in C^{r}(\mathbb{R})$, the following relation holds:
	\begin{equation}\label{dkf} 
		[D^{r}f](t)=\sum_{j=1}^{r}(-1)^{r-j}\frac{(2r-1-j)!}{2^{r-j}(r-j)!(j-1)!t^{2r-j}}f^{(j)}(t). 
	\end{equation}
\end{lemma}
We skip the proof of this  as well.

We can now compute $[D^{l-1}T_{l}](t)$ in terms of ordinary derivatives. 
\begin{lemma} For $l\geq 1$ and $f\in C^{l-1}(\mathbb{R})$, we have
	\begin{align}
		[D^{l-1}T_{l}](t)=(-1)^{l-1}\sum_{m=0}^{l-1}P_{m,l}(t)f^{(m)}(1-t)
	\end{align}
where
\begin{equation}\label{PMLL}
	P_{m,l}(t)=\sum_{n=m}^{l-1}E_{n,m,l}\frac{(1-t)^{n+2}}{t^{n+l-m}},\quad\text{with} \quad	E_{n,m,l}
	=\frac{(l-1)!}{2^{n-m}m!}\binom{l+1}{n+2}\binom{l-1+n-m}{n-m}.
\end{equation}
\end{lemma}
\begin{proof}
	Using Leibniz rule, Lemma \ref{Tlprelim} and Lemma \ref{fder}, we have
\begin{align}\label{Tl:derivative}
	&[D^{l-1}T_{l}](t)
	=D^{l-1}\left\{\frac{(1-t)^{l+1}}{t}\right\}f(1-t)+	\sum_{j=1}^{l-1}\binom{l-1}{j}D^{l-1-j}\left\{\frac{(1-t)^{l+1}}{t}\right\}D^{j}f(1-t)\nonumber\\
	&=D^{l-1}\left\{\frac{(1-t)^{l+1}}{t}\right\}f(1-t)+	\sum_{j=1}^{l-1}\binom{l-1}{j}\left[\sum_{m=1}^{j}(-1)^{j}\frac{(2j-1-m)!}{2^{j-m}(j-m)!(m-1)!}\frac{f^{(m)}(1-t)}{t^{2j-m}}\right]\nonumber\\
	&\times\left[(-1)^{l-1-j}(l-1-j)!\sum\limits_{s=0}^{l-1-j}\frac{1}{2^s} {l+1\choose l-1-j-s} {l-1-j+s\choose s}\frac{(1-t)^{j+s+2}}{t^{l-j+s}}\right]\nonumber\\
        &=D^{l-1}\left\{\frac{(1-t)^{l+1}}{t}\right\}f(1-t)+	(-1)^{l-1}\sum_{m=1}^{l-1}P_{m,l}(t)f^{(m)}(1-t),
\end{align}
where $P_{m,l}(t)$ is given by
\begin{align*}
	P_{m,l}(t)
	&=\sum_{j=m}^{l-1}\sum\limits_{s=0}^{l-1-j}\frac{(l-1)!}{j!2^{j+s-m}}\frac{(2j-1-m)!}{(j-m)!(m-1)!}{l+1\choose l-1-j-s} {l-1-j+s\choose s}\frac{(1-t)^{j+s+2}}{t^{l+j+s-m}}.
\end{align*}
Using the substitution $j+s=n$, the above expression becomes 
\begin{align*}
P_{m,l}(t)=\sum_{n=m}^{l-1}E_{n,m,l}\frac{(1-t)^{n+2}}{t^{n+l-m}},\quad m\geq 1,l\geq 2,
\end{align*}
where $E_{n,m,l}$ is given by
\Beq\label{enl:expression}
\begin{aligned} 
E_{n,m,l}&=\sum\limits_{p=m}^{n} \frac{(l-1)!(2p-1-m)!}{p! 2^{n-m}(p-m)!(m-1)!}{l+1\choose l-1-n}{l-1+n-2p\choose n-p}\\
&=\frac{(l-1)!}{2^{n-m}(m-1)!}{l+1\choose l-1-n}\sum\limits_{p=m}^{n}\frac{(2p-1-m)!}{p!(p-m)!}{l-1+n-2p\choose n-p}\\
&=\frac{(l-1)!}{2^{n-m}(m-1)!}{l+1\choose l-1-n}\sum\limits_{p=0}^{n-m}\frac{1}{2p+m}{2p+m\choose p}{l-1+n-2m-2p\choose n-m-p}. 
\end{aligned} 
\Eeq
Using Abel-Aigner identity, 
\[
\sum\limits_{k}\frac{r}{tk+r}{tk+r\choose k}{t(n-k)+s\choose n-k}={tn+r+s\choose n}, 
\]
we get that 
\[
\sum\limits_{p=0}^{n-m}\frac{1}{2p+m}{2p+m\choose p}{l-1+n-2m-2p\choose n-m-p}=\frac{1}{m}{n-m+l-1\choose n-m}.
\]
Therefore we get,
\Beq\label{Eqq4.6}
E_{n,m,l}=\frac{(l-1)!}{2^{n-m}m!}{l+1\choose l-1-n}{n-m+l-1\choose n-m}.
\Eeq
Thus $P_{m,l}(t)$ takes the form 
\begin{equation}\label{pml}
P_{m,l}(t)=\sum_{n=m}^{l-1}E_{n,m,l}\frac{(1-t)^{n+2}}{t^{n+l-m}},\quad\text{where} \quad	E_{n,m,l}
=\frac{(l-1)!}{2^{n-m}m!}\binom{l+1}{n+2}\binom{l-1+n-m}{n-m}.
\end{equation}
Using the above expression and 
Lemma \ref{Tlprelim}, it is easy to verify that
\begin{equation}\label{P_0l}
	P_{0,l}(t)=D^{l-1}\left\{\frac{(1-t)^{l+1}}{t}\right\}.
\end{equation}
Now we can immediately conclude the proof by combining  \eqref{Tl:derivative}, \eqref{pml}, and \eqref{P_0l}.
\end{proof}
The explicit expression for $\frac{\D}{\D t}[D^{r}G_{0,r}](t)$ in terms of the derivatives of $f$ can now be easily determined.
\begin{lemma}\label{lemmadDG0r} For any $r\geq 0$, we have
	\begin{align*}
		\left[\frac{\D}{\D t}D^{r}G_{0,r}\right](t)
		=2^{r}(-1)^{r}r!\sum_{m=0}^{r}\sum_{n=m}^{r}\sum_{l=n}^{r}\frac{1}{2^{n-m}m!}\binom{l+1}{n+1}\binom{l+n-m}{n-m}\frac{1}{t^{l}}\frac{(1-t)^{n+1}}{t^{n-m}}f^{(m)}(1-t).
	\end{align*}
\end{lemma}
\begin{proof}
From Lemma \ref{recurssion}, we have
 \begin{equation*}
	\left[\frac{\D}{\D t}D^{r}G_{0,r}\right](t)=2^{r}\left((-1)^{r}r!(1-t)f(1-t)+\sum_{l=1}^{r}(-1)^{r-l}\frac{r!}{l!}\left[\frac{\D}{\D t}D^{l-1}T_{l}\right](t)\right).
\end{equation*}
Now substituting the expression of $[D^{l-1}T_{l}](t)$, we get
\begin{align*}
&\frac{\D}{\D t}[D^{r}G_{0,r}](t)=2^{r}(-1)^{r}r!(1-t)f(1-t)+2^{r}(-1)^{r-1}\sum_{l=1}^{r}\sum_{m=0}^{l-1}\frac{r!}{l!}\frac{\D}{\D t}\left(	P_{m,l}(t)f^{(m)}(1-t)\right)\\
&=2^{r}(-1)^{r}r!\Bigg{\{}(1-t)f(1-t)-\sum_{l=1}^{r}\sum_{m=0}^{l-1}\frac{\lb P_{m,l}(t)\rb'}{l!}f^{(m)}(1-t)+\sum_{l=1}^{r}\sum_{m=1}^{l}\frac{P_{m-1,l}(t)}{l!}f^{(m)}(1-t)\Bigg{\}},
\end{align*}
where the last term follows reindexing the variable $m$. 
From \eqref{PMLL}
\begin{align*}
\lb P_{m,l}(t)\rb'=-\sum_{n=m}^{l-1}(n+2)E_{n,m,l}\frac{(1-t)^{n+1}}{t^{n+l-m}}-\sum_{n=m}^{l-1}(n+l-m)E_{n,m,l}\frac{(1-t)^{n+2}}{t^{n+l-m+1}}.
\end{align*}
Substituting the expressions of $P_{m,l}(t)$  and $\lb P_{m,l}(t)\rb'$ above, we have 
\begin{align}
&\frac{(-1)^r}{2^{r}r!}\frac{\D}{\D t}[D^{r}G_{0,r}](t)=(1-t)f(1-t)\nonumber\\
&+\sum_{l=1}^{r}\sum_{m=0}^{l-1}\frac{1}{l!}\left[\sum_{n=m}^{l-1}(n+2)E_{n,m,l}\frac{(1-t)^{n+1}}{t^{n+l-m}}+\sum_{n=m}^{l-1}(n+l-m)E_{n,m,l}\frac{(1-t)^{n+2}}{t^{n+l-m+1}}\right]f^{(m)}(1-t)\nonumber\\
&+\sum_{l=1}^{r}\sum_{m=1}^{l}\frac{1}{l!}\left[\sum_{n=m-1}^{l-1}E_{n,m-1,l}\frac{(1-t)^{n+2}}{t^{n+l-m+1}}\right]f^{(m)}(1-t)\nonumber\\
&=(1-t)f(1-t)\nonumber\\
&+\sum_{l=1}^{r}\sum_{m=0}^{l-1}\frac{1}{l!}\left[\sum_{n=m}^{l-1}(n+2)E_{n,m,l}+\sum_{n=m+1}^{l}(n-1+l-m)E_{n-1,m,l}\right]\frac{(1-t)^{n+1}}{t^{n+l-m}}f^{(m)}(1-t)\nonumber\\
&+\sum_{l=1}^{r}\sum_{m=1}^{l}\frac{1}{l!}\sum_{n=m}^{l}E_{n-1,m-1,l}\frac{(1-t)^{n+1}}{t^{n+l-m}}f^{(m)}(1-t).
\end{align}
Using \eqref{Eqq4.6}, we note that $E_{l,m,l}=E_{m-1,m,l}=E_{l,l,l}=E_{l-1,l,l}=0$, and
hence we get
\begin{align}\label{DrG-r:exxp}
    &\frac{(-1)^r}{2^{r}r!}\frac{\D}{\D t}[D^{r}G_{0,r}](t)=(1-t)f(1-t)\nonumber\\
&+\sum_{l=1}^{r}\sum_{m=0}^{l}\sum_{n=m}^{l}\frac{1}{l!}\left[(n+2)E_{n,m,l}+(n-1+l-m)E_{n-1,m,l}\right]\frac{(1-t)^{n+1}}{t^{n+l-m}}f^{(m)}(1-t)\nonumber\\
&+\sum_{l=1}^{r}\sum_{m=1}^{l}\sum_{n=m}^{l}\frac{1}{l!}E_{n-1,m-1,l}\frac{(1-t)^{n+1}}{t^{n+l-m}}f^{(m)}(1-t)\nonumber\\
&=\left[(1-t)+\sum_{l=1}^{r}\sum_{n=0}^{l}\frac{1}{l!}\left[(n+2)E_{n,0,l}+(n-1+l)E_{n-1,0,l}\right]\frac{(1-t)^{n+1}}{t^{n+l}}\right]f(1-t)\nonumber\\
&+\sum_{m=1}^{r}\sum_{l=m}^{r}\sum_{n=m}^{l}\frac{1}{l!}\left[(n+2)E_{n,m,l}+(n-1+l-m)E_{n-1,m,l}\right]\frac{(1-t)^{n+1}}{t^{n+l-m}}f^{(m)}(1-t)\nonumber\\
&+\sum_{m=1}^{r}\sum_{l=m}^{r}\sum_{n=m}^{l}\frac{1}{l!}E_{n-1,m-1,l}\frac{(1-t)^{n+1}}{t^{n+l-m}}f^{(m)}(1-t)\nonumber\\
&=\sum_{m=0}^{r}S_{m,r}(t)f^{(m)}(1-t),
\end{align}
where  $S_{m,r}(t)$ is given by
\begin{equation*}
   S_{0,r}(t)=(1-t)+\sum_{l=1}^{r}\sum_{n=0}^{l}\frac{1}{l!}\left[(n+2)E_{n,0,l}+(n+l-1)E_{n-1,0,l}\right]\frac{(1-t)^{n+1}}{t^{n+l}}
\end{equation*}
and for $1\leq m\leq r$, 
\begin{equation*}
    S_{m,r}(t)=\sum_{l=m}^{r}\sum_{n=m}^{l}\frac{1}{l!}\left[(n+2)E_{n,m,l}+(n-1+l-m)E_{n-1,m,l}+E_{n-1,m-1,l}\right]\frac{(1-t)^{n+1}}{t^{n+l-m}}.
\end{equation*}
 Now we will simplify the above expressions. Recall \eqref{PMLL}, the expression $S_{0,r}(t)$ takes the form
\begin{align}\label{coefff0}
S_{0,r}(t)
=&(1-t)+\sum_{l=1}^{r}\sum_{n=0}^{l}\frac{1}{2^n}\binom{l+1}{n+1}\frac{l+n}{l}\binom{l+n-1}{n}\frac{(1-t)^{n+1}}{t^{n+l}}\nonumber\\
=&(1-t)+\sum_{l=1}^{r}\sum_{n=0}^{l}\frac{1}{2^{n}}\binom{l+1}{n+1}\binom{l+n}{n}\frac{(1-t)^{n+1}}{t^{n+l}}\nonumber\\
=&\sum_{n=0}^{r}\sum_{l=n}^{r}\frac{1}{2^{n}}\binom{l+1}{n+1}\binom{l+n}{n}\frac{1}{t^{l}}\frac{(1-t)^{n+1}}{t^{n}}.
\end{align}
Recalling from \eqref{PMLL}, the expression $S_{m,r}(t)$ takes the form
\begin{align}\label{Smr}
    &S_{m,r}(t)\nonumber\\
    &=\sum_{l=m}^{r}\sum_{n=m}^{l}\frac{(l-1)!}{l!2^{n-m}m!}\binom{l+1}{n+1}\left[(n+2)\frac{l-n}{n+2}+2(n-m)+m\right]\binom{l-1+n-m}{n-m}\frac{(1-t)^{n+1}}{t^{n+l-m}}\nonumber\\
    &=\sum_{l=m}^{r}\sum_{n=m}^{l}\frac{1}{2^{n-m}m!}\binom{l+1}{n+1}\frac{l+n-m}{l}\binom{l-1+n-m}{n-m}\frac{(1-t)^{n+1}}{t^{n+l-m}}\nonumber\\
    &=\sum_{n=m}^{r}\sum_{l=n}^{r}\frac{1}{2^{n-m}m!}\binom{l+1}{n+1}\binom{l+n-m}{n-m}\frac{(1-t)^{n+1}}{t^{n+l-m}}.
\end{align}
Substituting \eqref{coefff0} and \eqref{Smr} in \eqref{DrG-r:exxp} gives the required expression for $\left[\frac{\D}{\D t}D^{r}G_{0,r}\right](t)$. This completes the proof.
\end{proof}
As a key component in proving our main theorem, we seek to express $G_{m,r}$ in terms of  $G_{m+j,0}$, $\left\lfloor \frac{r}{2}\right\rfloor\leq j\leq r$ and its lower order derivatives in the following lemma. We only require this lemma for the special case of $m=0$, that is, for $G_{0,r}$, however, we will prove a more general result below. 
\begin{lemma}\label{Gmrlemma}  Let $r\in\mathbb{N}_{0}$. Then 
\begin{equation}\label{conjec:main}
	\sum_{l=0}^{\left\lfloor \frac{r}{2}\right\rfloor}\frac{2^{2l}r!}{l!(r-2l)!}\frac{(m+r)!}{(m+r-l)!}[D^{r-2l}G_{m+r-l,0}](t)=4^{r}(m+1)\cdots(m+r)G_{m,r}(t)\quad \text{for all }m\in\mathbb{N}_{0},
\end{equation}
where $\left\lfloor \cdot\right\rfloor$ is the greatest integer function.
\end{lemma}
\begin{proof}
We will prove it by induction. Our induction argument is that, once we fix $r$ then \eqref{conjec:main} is true for all $m\geq 0$. The case  $r=0$  for any $m\geq 0$ is obvious.	For $r=1$ and any $m\geq 0$, the above relation gives $[DG_{m+1,0}](t)=4(m+1)G_{m,1}(t)$ which immediately follows from relation \eqref{rec:Gk0}. Now assume that it is true for all natural numbers up to $r$ and for any $m\geq 0$. Let us prove for $r+1$ and for any $m\geq 0$. Since \eqref{conjec:main} is true for $r$ and for any $m\geq 0$, we can replace $m$ by $m+1$
\begin{multline}\label{conjec}
    \sum_{l=0}^{\left\lfloor \frac{r}{2}\right\rfloor}\frac{2^{2l}r!}{l!(r-2l)!}\frac{(m+1+r)!}{(m+1+r-l)!}[D^{r-2l}G_{m+1+r-l,0}](t)=4^{r}(m+2)\cdots(m+1+r)G_{m+1,r}(t),
\end{multline}
Taking a $D$ derivative in \eqref{conjec}, we get
\begin{align*}
	&\sum_{l=0}^{\left\lfloor \frac{r}{2}\right\rfloor}\frac{2^{2l}r!}{l!(r-2l)!}\frac{(m+1+r)!}{(m+1+r-l)!}[D^{r+1-2l}G_{m+1+r-l,0}](t)\\
	&=4^{r}(m+2)\cdots(m+1+r)[DG_{m+1,r}](t)\\
	&=4^{r}(m+2)\cdots(m+1+r)\left[4(m+1)G_{m,r+1}(t)-2rG_{m+1,r-1}(t)\right]\quad (\text{using relation \eqref{rec:Gkj}})\\
	&=4^{r+1}(m+1)\cdots(m+1+r)G_{m,r+1}(t)-2r(m+1+r)4\left[4^{r-1}(m+2)\cdots(m+r)G_{m+1,r-1}(t)\right].
\end{align*}
Using the induction hypothesis, we get
\begin{multline*}
    	\sum_{l=0}^{\left\lfloor \frac{r}{2}\right\rfloor}\frac{2^{2l}r!}{l!(r-2l)!}\frac{(m+1+r)!}{(m+1+r-l)!}[D^{r+1-2l}G_{m+1+r-l,0}](t)=4^{r+1}(m+1)\cdots(m+1+r)G_{m,r+1}(t)
	-\\\sum_{l=0}^{\left\lfloor \frac{r-1}{2}\right\rfloor}\frac{2^{2l+3}r!}{l!(r-1-2l)!}\frac{(m+1+r)!}{(m+r-l)!}[D^{r-1-2l}h_{m+r-l}](t).
\end{multline*}
 Replacing  $l$ by $l-1$ in  summation on RHS
 \begin{multline*}
	\sum_{l=0}^{\left\lfloor \frac{r}{2}\right\rfloor}\frac{2^{2l}r!}{l!(r-2l)!}\frac{(m+r+1)!}{(m+r+1-l)!}[D^{r+1-2l}G_{m+r+1-l,0}](t)=4^{r+1}(m+1)\cdots(m+r+1)G_{m,r+1}
	-\\\sum_{l=1}^{\left\lfloor \frac{r+1}{2}\right\rfloor}\frac{2^{2l+1}r!}{(l-1)!(r+1-2l)!}\frac{(m+r+1)!}{(m+r+1-l)!}[D^{r+1-2l}h_{m+r+1-l}](t).
 \end{multline*}
 This gives
 \begin{equation*}
\sum_{l=0}^{\left\lfloor \frac{r+1}{2}\right\rfloor}\frac{2^{2l}(r+1)!}{l!(r+1-2l)!}\frac{(m+r+1)!}{(m+r+1-l)!}[D^{r+1-2l}G_{m+r+1-l,0}](t)=4^{r+1}(m+1)\cdots(m+r+1)G_{m,r+1},
 \end{equation*}
 which is exactly equal to \eqref{conjec:main} for $r+1$.
 This completes the proof.
\end{proof}
In the preceding lemma, 
the case corresponding to $m=0$ and $r=k$ is adequate for our purposes. That is,  
\begin{equation}\label{gokderiva}
	\sum_{l=0}^{\left\lfloor \frac{k}{2}\right\rfloor}\frac{2^{2l}k!}{(k-2l)!}\binom{k}{l}[D^{k-2l}G_{k-l,0}](t)=4^{k}k!G_{0,k}(t)\quad \text{for all } k\geq 0. 
	\end{equation}
Using the aforementioned lemma, we now prove in the following theorem that $\frac{\D}{\D t}[D^{2k}h_{k}](t)$ can be expressed solely in terms of the derivatives of $G_{0,j}(t)$'s. Next, we shall obtain our main theorem by incorporating the derivatives of $G_{0,j}(t)$.
\begin{theorem}\label{LemmaD2khk} For $k\geq 1$, we have
	\begin{equation}\label{g0kexp}
		\left[D^{2k}h_{k}\right](t)=\left[D^{2k}G_{k,0}\right](t)=k!4^{k}\sum_{j=0}^{k-1}(-1)^{j}\frac{(k-1+j)!}{(k-1-j)!j!}\left[D^{k-j}G_{0,k-j}\right](t).
	\end{equation}
\end{theorem}
\begin{proof} 
For $k=1$, 
	\begin{equation*}
		\left[D^{2}h_{1}\right](t)=4\left[DG_{0,1}\right](t).
	\end{equation*}
It is therefore valid for $k=1$. Now assume that it is true for any natural number up to $k$. We  now prove for $k+1$. We have to demonstrate 
		\begin{equation*}
		[D^{2k+2}G_{k+1,0}](t)=(k+1)!4^{k+1}\sum_{j=0}^{k}(-1)^{j}\frac{(k+j)!}{(k-j)!j!}[D^{k+1-j}G_{0,k+1-j}](t).
	\end{equation*}
	Recalling the relation \eqref{gokderiva} with $k$ replaced by $k+1$,
	\begin{equation*}
	\sum_{l=0}^{\left\lfloor \frac{k+1}{2}\right\rfloor}\frac{2^{2l}(k+1)!}{(k+1-2l)!}\binom{k+1}{l}[D^{k+1-2l}G_{k+1-l,0}](t)=4^{k+1}(k+1)!G_{0,k+1}(t).
	\end{equation*} 
Taking $(k+1)$ $D$-derivatives in the above equation, we obtain
		\begin{equation*}
		\sum_{l=0}^{\left\lfloor \frac{k+1}{2}\right\rfloor}\frac{2^{2l}(k+1)!}{(k+1-2l)!}\binom{k+1}{l}[D^{2(k+1-l)}G_{k+1-l,0}](t)=4^{k+1}(k+1)![D^{k+1}G_{0,k+1}](t).
	\end{equation*} 
	Now using the induction hypothesis for $[D^{2(k+1-l)}h_{k+1-l}](t)$, we get
	\begin{align*}
	&4^{k+1}(k+1)![D^{k+1}G_{0,k+1}](t)=[D^{2k+2}G_{k+1,0}](t)+\\	&	\sum_{l=1}^{\left\lfloor \frac{k+1}{2}\right\rfloor}\frac{2^{2l}(k+1)!}{(k+1-2l)!}\binom{k+1}{l}\left((k+1-l)!4^{k+1-l}\sum_{j=0}^{k-l}\frac{(-1)^{j}(k-l+j)!}{(k-l-j)!j!}[D^{k+1-l-j}G_{0,k+1-l-j}](t)\right)\\
    &=[D^{2k+2}G_{k+1,0}](t)+	((k+1)!)^2 4^{k+1}\sum_{l=1}^{\left\lfloor \frac{k+1}{2}\right\rfloor}\sum_{j=0}^{k-l}A_{j,l}[D^{k+1-l-j}G_{0,k+1-l-j}](t),
	\end{align*}
	where
	\begin{align*}
		A_{j,l}=(-1)^{j}\frac{(k-l+j)!}{(k+1-2l)!(k-l-j)!l!j!}.
	\end{align*}
	Using the substitution $l+j=m$, the above expression takes the form
	\begin{align}\label{bkm:eq}
	&[D^{2k+2}G_{k+1,0}](t)	\nonumber\\&=-\lb(k+1)!\rb^2 4^{k+1}\sum_{m=1}^{k}B_{k,m}[D^{k+1-m}G_{0,k+1-m}](t)+	4^{k+1}(k+1)![D^{k+1}G_{0,k+1}](t),
	\end{align}
where	the coefficients $B_{k,m}$ of $[D^{k+1-m}G_{0,k+1-m}](t)$ are given by 
	\begin{align*}
		B_{k,m}=\sum_{q=1}^{m}A_{m-q,q}&=\sum_{q=1}^{m}(-1)^{m-q}\frac{(k-2q+m)!}{(k+1-2q)!(k-m)!q!(m-q)!}\\
		&=\frac{(-1)^{m}}{(k-m)!m!}\sum_{q=1}^{m}(-1)^{q}\frac{(k-2q+m)!}{(k+1-2q)!}\binom{m}{q}\\
		&=\frac{(-1)^{m}}{(k-m)!m!}\sum_{q=1}^{m}(-1)^{q}P(q)\binom{m}{q},
	\end{align*}
where $P(q)=\frac{(k-2q+m)!}{(k+1-2q)!}=(k-2q+2)\cdots(k-2q+m)$ is a polynomial of degree $m-2$ in $q$. Now using the fact that
\begin{equation*}
	\sum_{j=0}^{n}(-1)^{j}P(j)\binom{n}{j}=0,
\end{equation*}
for any polynomial $P$ of degree less than $n$, we deduce that
\begin{equation*}
B_{k,m}=\frac{(-1)^{m}}{(k-m)!m!}\left[-P(0)\binom{m}{0}\right]=-\frac{(-1)^{m}(k+m)!}{(k-m)!m!(k+1)!}.
\end{equation*}
Substituting $B_{k,m}$ back into the equation \eqref{bkm:eq}, we get
	\begin{align*}
	&[D^{2k+2}G_{k+1,0}](t)\\
    &=4^{k+1}(k+1)!\sum_{m=1}^{k}\frac{(-1)^{m}(k+m)!}{(k-m)!m!}[D^{k+1-m}G_{0,k+1-m}](t)+4^{k+1}(k+1)![D^{k+1}G_{0,k+1}](t)\\
	&=4^{k+1}(k+1)!\sum_{m=0}^{k}\frac{(-1)^{m}(k+m)!}{(k-m)!m!}[D^{k+1-m}G_{0,k+1-m}](t).
\end{align*}
Taking a $\frac{\D}{\D t}$ derivative, we immediately obtain \eqref{g0kexp} for $k+1$.
This completes the proof.
\end{proof}
\section{Proofs}\label{sec:proof of main results}
We are now ready with all the tools to prove our first main result.
\begin{proof}[Proof of Theorem \ref{mainthm1}]
	For $k=0$, it is obvious. For $k\geq 1$,
from Theorem \ref{LemmaD2khk}, we have
\begin{equation*}
		\left[\frac{\D}{\D t}D^{2k}h_{k}\right](t)=k!4^{k}\sum_{j=0}^{k-1}(-1)^{j}\frac{(k-1+j)!}{(k-1-j)!j!}\left[\frac{\D}{\D t}D^{k-j}G_{0,k-j}\right](t).
\end{equation*}
Recalling the relation for $\left[\frac{\D}{\D t}D^{k-j}G_{0,k-j}\right](t)$ from Lemma \ref{lemmadDG0r}, the above expression becomes
\begin{align*}
&\frac{1}{(-1)^{k}4^{k}k!}\left[\frac{\D}{\D t}D^{2k}h_{k}\right](t)\nonumber\\
			&=\sum_{j=0}^{k-1}\frac{2^{k-j}(k-1+j)!(k-j)!}{(k-1-j)!j!}\sum_{m=0}^{k-j}\sum_{n=m}^{k-j}\sum_{l=n}^{k-j}\frac{1}{2^{n-m}m!}\binom{l+1}{n+1}\binom{l+n-m}{n-m}\frac{1}{t^{l}}\frac{(1-t)^{n+1}}{t^{n-m}}f^{(m)}(1-t)\nonumber\\
            &=\sum_{j=0}^{k}\sum_{m=0}^{k-j}\sum_{n=m}^{k-j}\sum_{l=n}^{k-j}\frac{2^{k-j}(2j)!(k-j)!}{2^{n-m}m!j!}\binom{k-1+j}{k-1-j}\binom{l+1}{n+1}\binom{l+n-m}{n-m}\frac{1}{t^{l}}\frac{(1-t)^{n+1}}{t^{n-m}}f^{(m)}(1-t),
\end{align*}
where the convention from \cite{Egorychev} using contour integration for $\binom{k-1+j}{k-1-j}$ allows to extend the sum up to $k$ in the $j$ variable. By changing the order of summation several times, we have 
\begin{align}\label{mainest}
    &\frac{1}{(-1)^{k}4^{k}k!}\left[\frac{\D}{\D t}D^{2k}h_{k}\right](t)\nonumber\\   &=\sum_{m=0}^{k}\sum_{n=m}^{k}\sum_{l=n}^{k}\sum_{j=0}^{k-l}\frac{2^{k-j}(2j)!(k-j)!}{2^{n-m}m!j!}\binom{k-1+j}{k-1-j}\binom{l+1}{n+1}\binom{l+n-m}{n-m}\frac{1}{t^{l}}\frac{(1-t)^{n+1}}{t^{n-m}}f^{(m)}(1-t).
\end{align}
It is straightforward to see that the inner sum: 
\begin{align}\label{clainkl}
	\sum_{j=0}^{k-l}2^{k-j}\frac{(2j)!(k-j)!}{j!}\binom{k-1+j}{k-1-j}=\frac{2^{l}(2k-l)!}{(k-l)!}.
\end{align}
Combining the equations \eqref{mainest} and \eqref{clainkl}, we get
\begin{multline*}
	\left[\frac{\D}{\D t}D^{2k}h_{k}\right](t)\\
	=(-1)^{k}k!4^{k}\sum_{m=0}^{k}\sum_{n=m}^{k}\sum_{l=n}^{k}2^{l+m-n}\frac{(2k-l)!}{(k-l)!m!}\binom{l+1}{n+1}\binom{l+n-m}{n-m}\frac{1}{t^{l}}\frac{(1-t)^{n+1}}{t^{n-m}}f^{(m)}(1-t).
\end{multline*}
This completes the proof.
\end{proof}
Using the explicit inversion for radial functions, we can now prove our second main theorem.
In the following proof, we use the computation and methodology employed in \cite{Salman_Article,SRTrange}.
\begin{proof}[Proof of Theorem \ref{mainthm2}]
The spherical harmonics expansion of $g=\mathcal{R}f$ is given by    
\begin{equation*}
    g(\theta,t)=\sum_{q=0}^{\infty}\sum_{s=1}^{d_q} g_{q,s}(t)Y_{q,s}(\theta),
\end{equation*}
where 
\begin{equation}\label{gqs}
    g_{q,s}(t)=\frac{\omega_{n-2}}{t^{n-2}\omega_{n-1}C_q^{\frac{n-2}{2}}(1)}\int_{1-t}^{1}u^{n-2}f_{q,s}(u)C_{q}^{\frac{n-2}{2}}\left(\frac{1+u^2-t^2}{2u}\right)\left(1-\left(\frac{1+u^2-t^2}{2u}\right)^{2}\right)^{\frac{n-3}{2}}du,
\end{equation}
where $C_q^{\frac{n-2}{2}}(t)$ are the Gegenbauer polynomials; see \cite{Reimer2003}, with  $C_q^{\frac{n-2}{2}}(1)=\dfrac{\Gamma(n-2+q)}{\Gamma(n-2)q!}$  and $f_{q,s}(r)$ is given in \eqref{dqandfqs}. Using the following relation established in \cite[Equation 3.31]{SRTrange},
\begin{equation*}
   C_{q}^{\frac{n-2}{2}}\left(\frac{1+u^2-t^2}{2u}\right)=K_{q,n}\left(1-\left(\frac{1+u^2-t^2}{2 u}\right)^2\right)^{-\frac{n-3}{2}}(-u)^q D^q\left(1-\left(\frac{1+u^2-t^2}{2 u}\right)^2\right)^{q+\frac{n-3}{2}},
\end{equation*}
where $K_{q,n}=\dfrac{(-1)^q \Gamma\left(\frac{n-1}{2}\right) \Gamma(n-2+q)}{2^q q!\Gamma( n-2) \Gamma\left(q+\frac{n-1}{2}\right)}$, the expression \eqref{gqs} takes the form
\begin{align*}
    g_{q, s}(t)&=\frac{K_{q,n}(-1)^q \omega_{n-2}}{t^{n-2}\omega_{n-1} C_{q}^{\frac{n-2}{2}}(1)} \int_{1-t}^1 u^{q+n-2} f_{q, s}(u) D^q\left(1-\frac{\left(1+u^2-t^2\right)^2}{4 u^2}\right)^{q+\frac{n-3}{2}} \mathrm{~d} u\\
    &=\frac{\left(\frac{n-3}{2}\right)! \omega_{n-2}}{2^{q}4^{q+\frac{n-3}{2}}\left(q+\frac{n-3}{2}\right)!t^{n-2}\omega_{n-1} } D^{q}\int_{1-t}^1 u^{1-q} f_{q, s}(u)\left(4u^2-(1+u^2-t^2)^2\right)^{q+\frac{n-3}{2}} \mathrm{~d} u.
\end{align*}
Going back to the definition of $Q(t,u)$, noting that $k:=\frac{n-3}{2}$ and
 denoting $\widetilde{f}_{q,s}(u)=\dfrac{f_{q,s}(u)}{u^{q}}$, the above expression becomes
\begin{align*}
    \dfrac{\omega_{2k+2}t^{2k+1}}{\omega_{2k+1}}g_{q, s}(t)
    &=\frac{k!}{ 2^{q}4^{q+k}(q+k)!} D^{q}\int_{1-t}^1 u \widetilde{f}_{q, s}(u)[Q(t,u)]^{q+k} \mathrm{~d} u.
\end{align*}
If we set 
\begin{equation*}
    h_{q,s}(t):=\dfrac{\omega_{2k+2}t^{2k+1}}{\omega_{2k+1}}g_{q, s}(t)\quad \text{and}\quad \phi_{q,s}(t):=\int_{1-t}^1 u \widetilde{f}_{q, s}(u)[Q(t,u)]^{q+k} \mathrm{~d} u,
\end{equation*}
then we have
\begin{equation}\label{aphirel}
    h_{q,s}(t)=\frac{k!}{ 2^{q}4^{q+k}(q+k)!} [D^{q}\phi_{q,s}](t).
\end{equation}
Given that $\widetilde{f}_{q, s}$ is radial, we may use Theorem \ref{mainthm1} to obtain
\begin{multline*}
	\frac{1}{(-1)^{q+k}(q+k)!4^{q+k}}\left[\frac{\D}{\D t}D^{2q+2k}\phi_{q,s}\right](t)\\
=\sum_{m=0}^{q+k}\sum_{n=m}^{q+k}\sum_{l=n}^{q+k}2^{l+m-n}\frac{(2q+2k-l)!}{(q+k-l)!m!}\binom{l+1}{n+1}\binom{l+n-m}{n-m}\frac{1}{t^{l}}\frac{(1-t)^{n+1}}{t^{n-m}}\widetilde{f}_{q,s}^{(m)}(1-t).
\end{multline*}
Using the relation \eqref{aphirel}, we get
\begin{multline*}
	\frac{2^{q}}{(-1)^{q+k}k!}\left[\frac{\D}{\D t}D^{q+2k}h_{q,s}\right](t)\\
=\sum_{m=0}^{q+k}\sum_{n=m}^{q+k}\sum_{l=n}^{q+k}2^{l+m-n}\frac{(2q+2k-l)!}{(q+k-l)!m!}\binom{l+1}{n+1}\binom{l+n-m}{n-m}\frac{1}{t^{l}}\frac{(1-t)^{n+1}}{t^{n-m}}\widetilde{f}_{q,s}^{(m)}(1-t).
\end{multline*}
This completes the proof.
\end{proof}
\section{Analytical Solutions}\label{sec:analytical sol}
We provide the analytical solution to the ODE obtained from Theorem \ref{mainthm1} for a few specific cases in this section.
\begin{example}[Dimension $n=3$] For $n=3$ i.e., $k=0$, Theorem \ref{mainthm1} gives
\begin{equation*}
    \frac{\D}{\D t} h_{0}(t)=(1-t)f(1-t)\implies f(t)=\frac{h_{0}^{\prime}(1-t)}{t}.
\end{equation*}   
\end{example}
\begin{example}[Dimension $n=5$] For $n=5$ i.e., $k=1$, Theorem \ref{mainthm1} gives
 \begin{align*}
	\left[\frac{\D}{\D t}D^{2}h_{1}\right](t)
=-8\left(\frac{(1-t)^{2}}{t}f^{(1)}(1-t)+\frac{(1-t)(1+t+t^{2})}{t^2}f(1-t)\right),
\end{align*}
and the solution of this ODE is given by 
\begin{equation*}
	f(1-t)=\frac{t}{8(1-t)^3}e^{-t}\int_{\epsilon}^{t}(1-u)e^{u}\left[\frac{\D}{\D u}D^{2}h_{1}\right](u)\D u.
\end{equation*}
\end{example}
\begin{example}[Dimension $n=7$] For $n=7$ i.e., $k=2$, Theorem \ref{mainthm1} gives
    \begin{equation*}\label{n=7:homode}
    \left[\frac{\D}{\D t}D^{4}h_{2}\right](t)=
	128\left(\frac{(1-t)^{3}}{t^2}f^{(2)}(1-t)+3\frac{(1-t)^2(1+t+t^2)}{t^3}f^{(1)}(1-t)+3\frac{(1-t^{5})}{t^4}f(1-t)\right).
\end{equation*}
The straightforward computation gives the solution of the above ODE as:
\begin{multline*}
     f(1-t)=-f_{1}(t)\int_{\epsilon}^{t}\frac{f_{2}(s)}{\mathrm{W}(f_{1},f_{2})}\frac{s^2}{128(1-s)^3}\left[\frac{\D}{\D s}D^{4}h_{2}\right](s)\D s+\\f_{2}(t)\int_{\epsilon}^{t}\frac{f_{1}(s)}{\mathrm{W}(f_{1},f_{2})}\frac{s^2}{128(1-s)^3}\left[\frac{\D}{\D t}D^{4}h_{2}\right](s)\D s,
\end{multline*}
where $f_{1}(t)$ and $f_{2}(t)$ are complementary solutions of homogeneous ODE and given by  $ f_{1}(t)=e^{-\frac{3}{2}t}\frac{\left[(t-2t^2)\cos\left(\frac{\sqrt{3}}{2}t\right)+\sqrt{3}t\sin\left(\frac{\sqrt{3}}{2}t\right)\right]}{(1-t)^{5}}$ and $f_{2}(t)=e^{-\frac{3}{2}t}\frac{\left[\sqrt{3}t\cos\left(\frac{\sqrt{3}}{2}t\right)-(t-2t^2)\sin\left(\frac{\sqrt{3}}{2}t\right)\right]}{(1-t)^{5}}$, and $\mathrm{W}(f_{1},f_{2})$ is the Wronskian given by
$\mathrm{W}(f_{1},f_{2})=\frac{2\sqrt{3}t^3e^{-\frac{3}{2}t}}{(1-t)^{9}}$.
\end{example}
\begin{remark}
To determine the complementary (homogeneous) solution of the ordinary differential equation (ODE) in arbitrary odd dimensions $n$, one can follow the method in \cite[Section 15.5]{Ince}  and then using the Duhamel's principle to obtain the desired solution; see \cite[Chapter 3]{Coddington_Book}.
\end{remark}
\section{Numerical Simulations}\label{sec:num simulations}
\begin{figure}[H]
    \centering
       \begin{subfigure}[b]{0.35\textwidth}
        \centering
    \includegraphics[width=\textwidth]{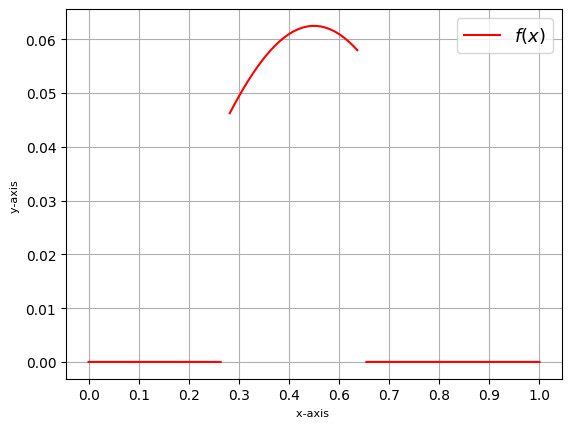}  
    \end{subfigure}
    \begin{subfigure}[b]{0.35\textwidth}
        \centering
        \includegraphics[width=\textwidth]{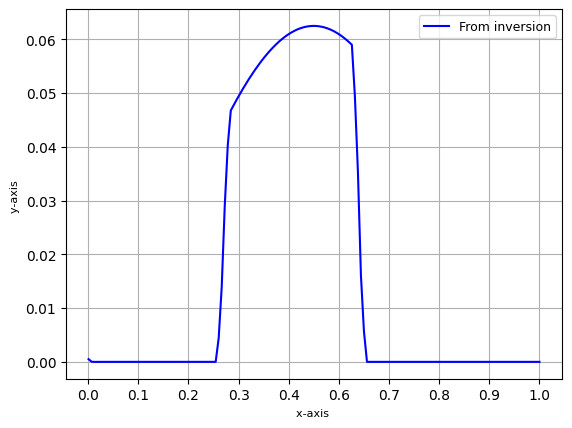}
    \end{subfigure}
    \caption{Dimension $n=3$, $f(x)=
\begin{cases}
     x^2(1-x)^2 &\text{if } 0.3<x<0.6,\\
     0 & \text{else}
 \end{cases}    $ on $[0.001,1]$, 100 node points. Here again, there is no ODE, but we consider a discontinuous  radial function.}
 \label{fig:1}
\end{figure}

\begin{figure}[H]
    \centering
       \begin{subfigure}[b]{0.35\textwidth}
        \centering
        \includegraphics[width=\textwidth]{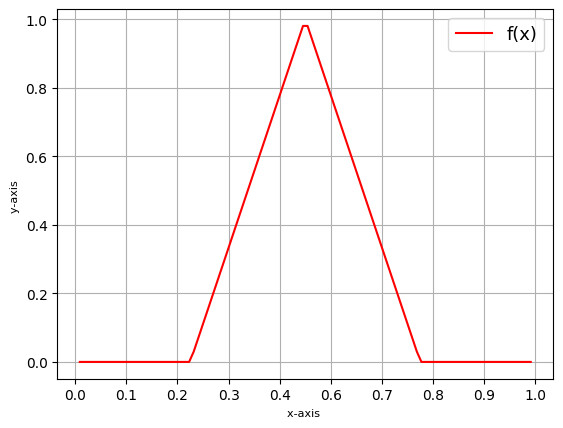} 
    \end{subfigure}
    \begin{subfigure}[b]{0.35\textwidth}
        \centering
        \includegraphics[width=\textwidth]{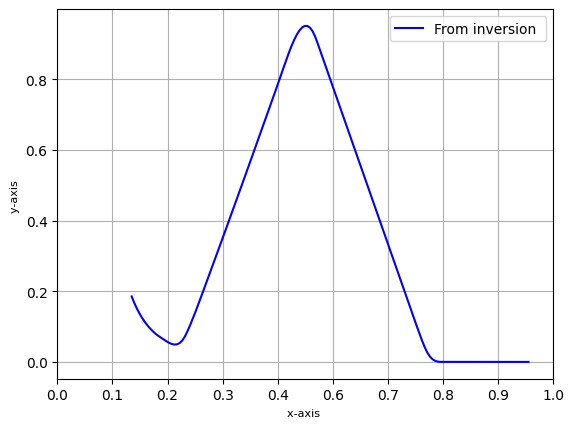}
    \end{subfigure}
    \caption{Dimension $n=5$, $f(x)=
\begin{cases}
     4x-1 &\text{if } 1/4<x<1/2,\\
     3-4x &\text{if } 1/2\leq x <3/4,\\
     0 &\text{else}
 \end{cases} $ on $[0.15,.95]$,  100 node points. Due to the singularity of the ODE at the origin and the resulting lack of smoothness of the function, the reconstruction is inaccurate near the origin. However, we are able to reconstruct the function over more than 70\%  of its domain.}
 \label{fig:2}
\end{figure}

\begin{figure}[H]
    \centering
       \begin{subfigure}[b]{0.35\textwidth}
        \centering   
        \includegraphics[width=\textwidth]{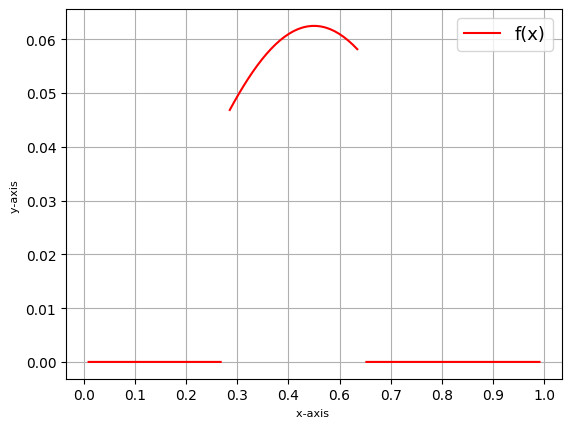}      
    \end{subfigure}
    \begin{subfigure}[b]{0.35\textwidth}
        \centering
        \includegraphics[width=\textwidth]{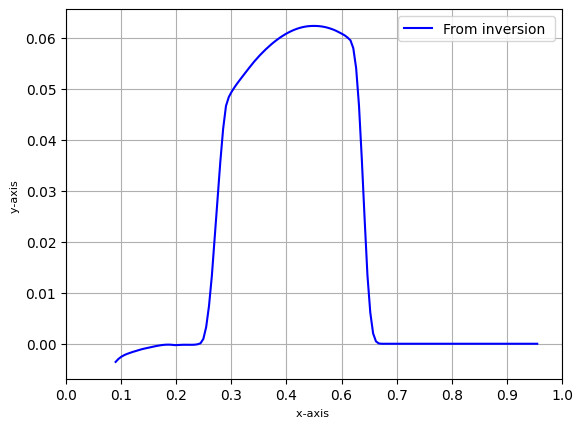}
 \end{subfigure}
   \caption{Dimension $n=5$, $f(x)=
\begin{cases}
     x^2(1-x)^2 &\text{if } 0.3<x<0.6\\
     0 & \text{else}
 \end{cases} $ on $[0.1,.95]$,  100 node points. Due to the singularity of the ODE at the origin and the resulting lack of smoothness of the function, the reconstruction is inaccurate near the origin. However, we are able to reconstruct the function over more than 80\% of its domain.}
 \label{fig:3}
\end{figure}

\begin{figure}[H]
    \centering
       \begin{subfigure}[b]{0.35\textwidth}
        \centering
        \includegraphics[width=\textwidth]{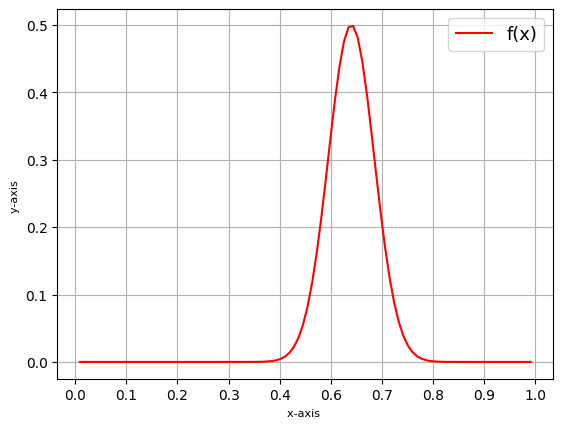} 
    \end{subfigure}
    \begin{subfigure}[b]{0.35\textwidth}
        \centering
        \includegraphics[width=\textwidth]{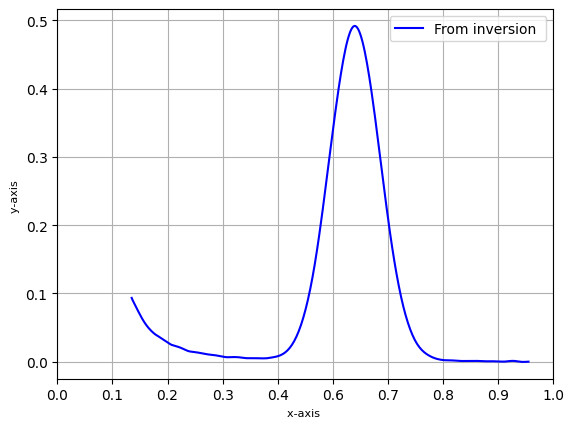}
    \end{subfigure}
\caption{Dimension $n=5$, $f(x)=\frac{1}{2}e^{-\frac{{(x-.6)^2}}{2(.05)^2}}$ with uniform noise ranges between $(10^{-7},-10^{-7})$ on $[0.15,.95]$, 100 node points. The presence of noise and a singularity in the ODE makes the reconstruction of the  function near the origin far from accurate.}
 \label{fig:9}
\end{figure}

This section presents the numerical simulations for the explicit inversion of SMT for both general and radial functions to illustrate the efficiency of the proposed inversion formula
in order to support our theoretical results. The python codes for the following simulations are available on \href{https://github.com/Pradipta97/SRTPartialData.git}{Github}. We remark that the inversion algorithm applies for SMT of compactly supported smooth functions in the unit ball. However in numerical simulations, we tested the reconstruction algorithm on SMT of non-differentiable functions with finitely many discontinuities and with noise added as well; see figures \ref{fig:1}, \ref{fig:2}, \ref{fig:3} and \ref{fig:9}.
We only considered low noise levels and addressing high noise levels requires a more thorough investigation. 

For the reconstruction of a general function from its SMT, the algorithm is as follows:
    \begin{enumerate}
        \item Given the data $g(\theta,t)$, consider its spherical harmonic expansion, 
        $$
g(\theta,t)=\sum_{q=0}^{\infty}\sum_{s=1}^{d_q} g_{q,s}(t)Y_{q,s}(\theta).
$$
\item Truncate the infinite series and compute  $g_{q,s}(t)$ using  \eqref{dqandfqs}.  
\item Compute $h_{q,s}(t)$ by the formula $h_{q,s}(t):=\dfrac{\omega_{2k+2}t^{2k+1}}{\omega_{2k+1}}g_{q, s}(t)\quad$ and use the ODE in Theorem \ref{mainthm2} to recover $\widetilde{f}_{q,s}(t)$.
\item We get $f_{q,s}(t)$ using $\widetilde{f}_{q,s}(t)=\dfrac{f_{q,s}(t)}{t^{q}}$. Now using the truncated spherical harmonics expansion of $f(t)$, we recover the approximate $f$.
    \end{enumerate}
\begin{figure}[H]
    \centering
       \begin{subfigure}[b]{0.35\textwidth}
        \centering
        \includegraphics[width=\textwidth]{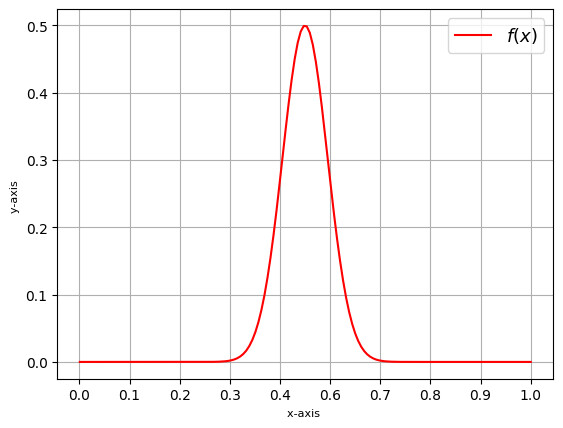} 
    \end{subfigure}
    \begin{subfigure}[b]{0.35\textwidth}
        \centering
        \includegraphics[width=\textwidth]{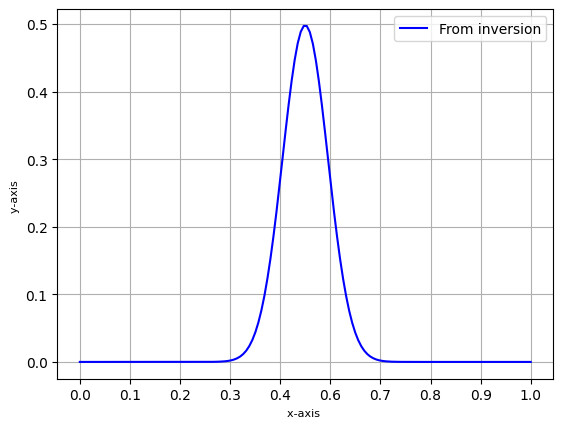}  
    \end{subfigure}
    \caption{Dimension $n=3$, $f(x)=\frac{1}{2}e^{-\frac{{(x-.5)^2}}{2(.05)^2}}$ on $[0.0001,1]$, 150 node points. In this case, there is no associated ODE; instead, the solution is given by an explicit formula. This enables us to reconstruct the function in a neighborhood arbitrarily close to the origin.}
\end{figure}
  \begin{figure}[H]
    \centering
       \begin{subfigure}[b]{.35\textwidth}
        \centering
        \includegraphics[width=\textwidth]{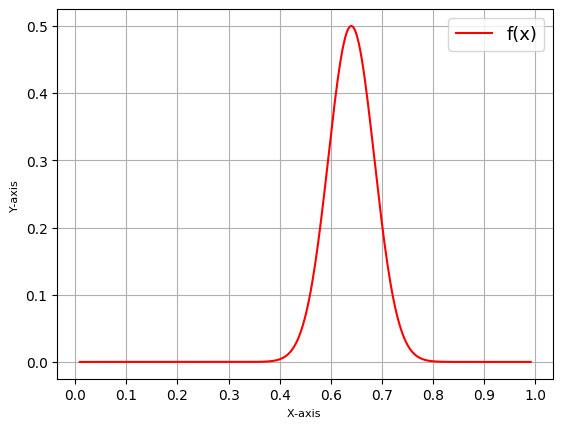} 
    \end{subfigure}
    \begin{subfigure}[b]{0.35\textwidth}
        \centering
        \includegraphics[width=\textwidth]{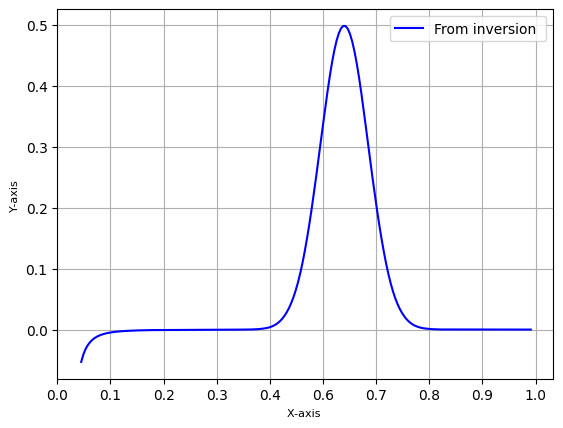}
    \end{subfigure}
    \caption{Dimension $n=5$, $f(x)=\frac{1}{2}e^{-\frac{{(x-.6)^2}}{2(.05)^2}}$ on $[0.05,0.99]$, 300 node points. Due to the singularity of the ODE at the origin, reconstruction of the function near the origin is inaccurate. }
\end{figure} 
\begin{figure}[H]
    \centering
       \begin{subfigure}[b]{0.35\textwidth}
        \centering
        \includegraphics[width=\textwidth]{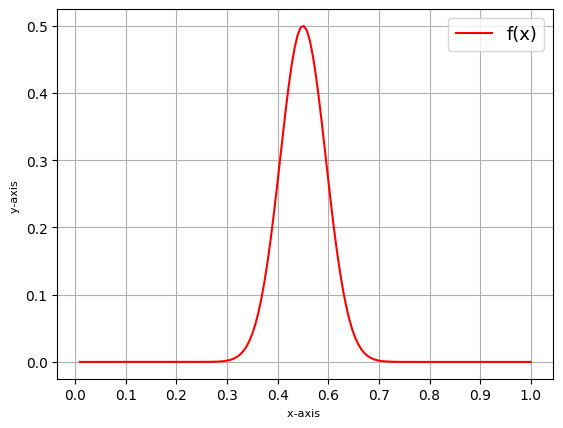} 
    \end{subfigure}
    \begin{subfigure}[b]{0.35\textwidth}
        \centering
        \includegraphics[width=\textwidth]{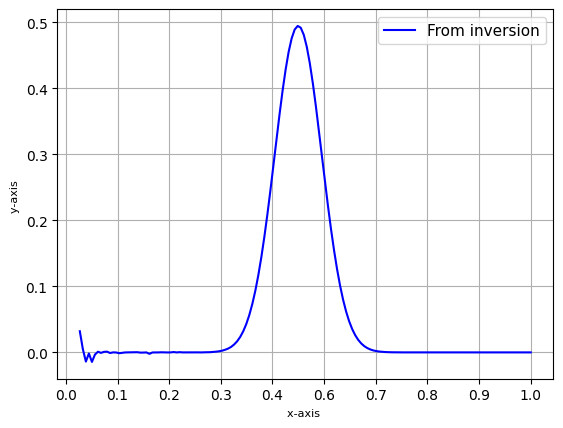}
    \end{subfigure}   
\caption{Dimension $n=5$, $f(x)=\frac{1}{2}e^{-\frac{{(x-.5)^2}}{2(.05)^2}}$ obtained from explicit solution using  Duhamel's principle on $[0.03,1]$, 150 node points.}
\end{figure}   
\begin{figure}[H]
    \centering
       \begin{subfigure}[b]{0.35\textwidth}
        \centering
        \includegraphics[width=\textwidth]{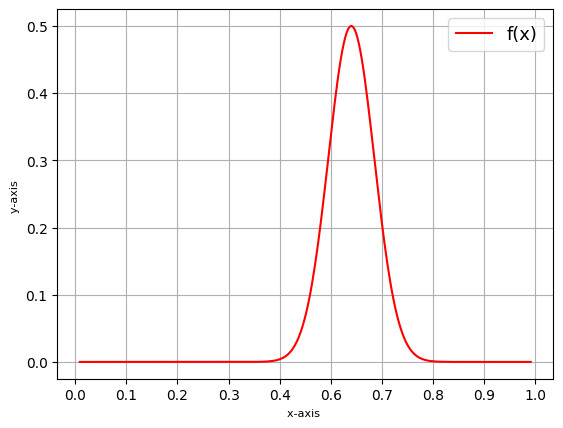} 
    \end{subfigure}
    \begin{subfigure}[b]{0.35\textwidth}
        \centering
        \includegraphics[width=\textwidth]{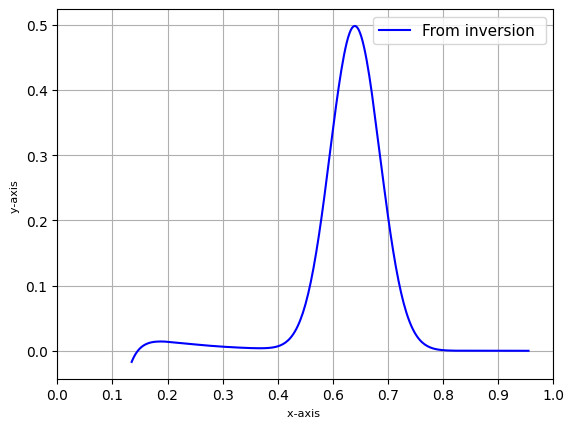}
    \end{subfigure}
\caption{Dimension $n=7$, $f(x)=\frac{1}{2}e^{-\frac{{(x-.6)^2}}{2(.05)^2}}$ on $[0.15,.95]$, 300 node points. Here, we solve a second-order ODE that is singular at the origin. Since computing higher-order derivatives introduces significant round-off errors, we are able to reconstruct the function over only 80\% of the domain.}
\end{figure}   
\begin{figure}[H]
    \centering
       \begin{subfigure}[b]{0.35\textwidth}
        \centering
        \includegraphics[width=\textwidth]{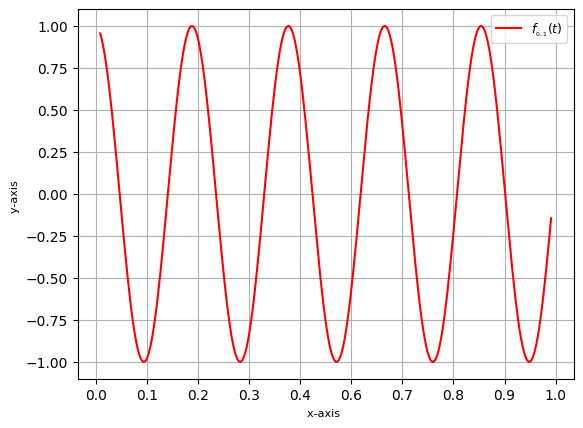} 
    \end{subfigure}
    \begin{subfigure}[b]{0.35\textwidth}
        \centering
        \includegraphics[width=\textwidth]{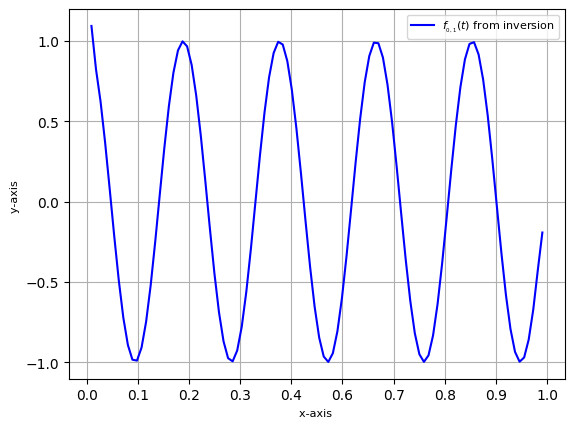}
    \end{subfigure}
    \begin{subfigure}[b]{0.35\textwidth}
        \centering
        \includegraphics[width=\textwidth]{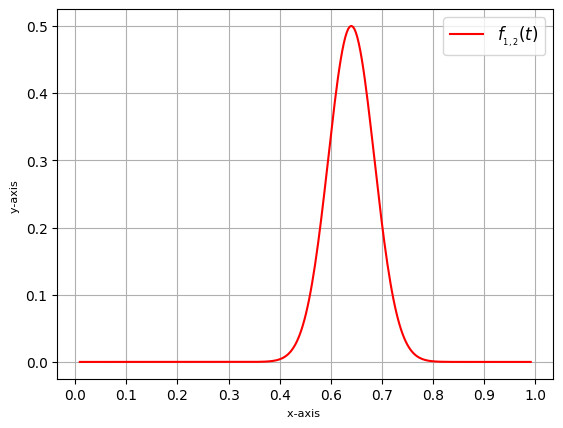}
            \end{subfigure}
    \begin{subfigure}[b]{0.35\textwidth}
        \centering
        \includegraphics[width=\textwidth]{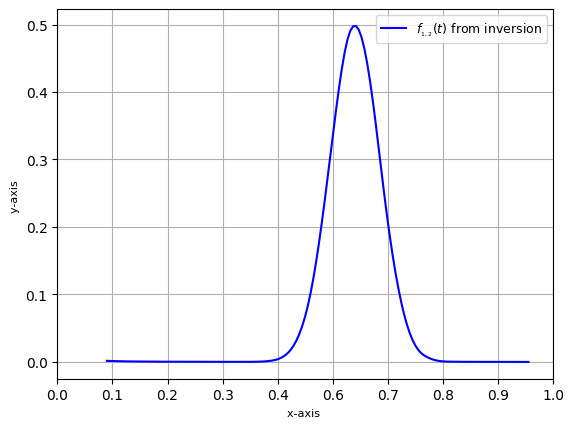}
    \end{subfigure}
 \caption{General function for dimension n=3, $f(x)=f_{0,1}(t)Y_{0,1}(\theta)+f_{1,2}(t)Y_{1,2}(\theta)$, where $x=t\theta$  $f_{0,1}(t)=\cos(50t)$, $f_{1,2}(t)=\frac{1}{2}e^{-\frac{{(x-.7)^2}}{2(.05)^2}}$ on $[0.01,.99]$ and $[0.1,0.95]$, 300 node points.}
    \label{fig:8}
\end{figure}   
\section{Acknowledgements}
PC acknowledges the support of the Department of Atomic Energy,  Government of India, for PhD Fellowship.

VPK would like to thank the Isaac Newton Institute for Mathematical Sciences, Cambridge, UK, for support and hospitality during \emph{Rich and Nonlinear Tomography - a multidisciplinary approach} in 2023 (supported by EPSRC Grant Number EP/R014604/1).
Additionally, VPK acknowledges the support of the Department of Atomic Energy,  Government of India, under
Project No.  12-R\&D-TFR-5.01-0520.

AT acknowledges the supported by the Anusandhan National Research Foundation (ANRF), under the scheme National
Post-Doctoral Fellowship, file no. PDF/2025/005843.

\bibliographystyle{plain}
%\bibliography{references.bib}
\end{document}